\documentclass[10pt,leqno]{article}

\usepackage{geometry}
\usepackage{amsthm,authblk,amsmath}
\usepackage[T1]{fontenc}
\usepackage{amsfonts}
\usepackage{graphicx}
\usepackage{epstopdf}
\usepackage{enumitem}
\usepackage{algorithmic}
\usepackage{amssymb}
\usepackage{dsfont}
\usepackage{subcaption}
\ifpdf
  \DeclareGraphicsExtensions{.eps,.pdf,.png,.jpg}
\else
  \DeclareGraphicsExtensions{.eps}
\fi

\theoremstyle{theorem}
\newtheorem{theorem}{Theorem}

\newtheorem{conjecture}[theorem]{Conjecture}
\newtheorem{definition}[theorem]{Definition}

\usepackage{xcolor}

\usepackage{hyperref}
\hypersetup{colorlinks=true, breaklinks=true, urlcolor= blue, linkcolor= teal, citecolor=teal,
  pdftitle={Spectral Gaps on Large Hyperbolic Surfaces}, 
  pdfauthor={Laura Monk and Frédéric Naud}, 
  pdfsubject={}}

\usepackage{amsopn}
\renewcommand{\Re}{{\rm Re}}
\renewcommand{\Im}{{\rm Im}}
\newcommand{\R}{\mathbb{R}}
\newcommand{\C}{\mathbb{C}}
\newcommand{\Z}{\mathbb{Z}}
\newcommand{\N}{\mathbb{N}}

\renewcommand{\H}{\mathbb{H}^2}
\renewcommand{\d}{\; \text{d}}
\newcommand{\dvwp}{\d \mathrm{Vol}_g^{\mathrm{WP}}}
\newcommand{\pwp}{\mathbb{P}_g^{\mathrm{WP}}}

\usepackage[capitalize]{cleveref}

\begin{document}

\title{\Large Spectral Gaps on Large Hyperbolic Surfaces}

\author[$\dagger$]{Laura Monk}
\author[$\star$]{Frédéric Naud}

\affil[$\dagger$]{University of Bristol, United Kingdom.}
\affil[$\star$]{Institut Mathématique de Jussieu, Université Pierre et Marie Curie, 4 place Jussieu,
    75252 Paris Cedex 05, France.}

\date{}

\maketitle

\begin{abstract} In this expository paper, we review the history and the recent breakthroughs in the
  spectral theory of large volume hyperbolic surfaces. More precisely, we focus mostly on the
  investigation of the first non-trivial eigenvalue $\lambda_1$ and its possible behaviour in the
  large volume regime.
\end{abstract}

\section{Introduction: from Selberg to Mirzakhani, the spectral gap of hyperbolic surfaces.}

\label{sec:intro}

\subsection{Spectral gap: history, expected behaviour, interactions.}

\subsubsection{Definition.}

In what follows, a {\it hyperbolic surface} $X$ will be a quotient $X=\Gamma \backslash \H$ of the
hyperbolic plane $\H$, endowed with the standard metric of Gauss curvature $-1$, by a Fuchsian group
of isometries $\Gamma$. The positive Laplace operator $\Delta_X$ acting as an unbounded operator on
$L^2$ is essentially self-adjoint. If $\Gamma$ is {\it co-compact or co-finite}, the bottom of the
spectrum of $\Delta_X$ is $0$ (which corresponds to the constant eigenfunctions) and the spectrum
below $1/4$ consists of finitely many eigenvalues, with finite multiplicity. Here the magical number
$1/4$ doesn't show up by accident: it is the bottom of the spectrum of the Laplacian $\Delta_{\H}$
acting on $L^2(\H)$, endowed with the hyperbolic volume \cite{mckean1970}. In this paper, we will
denote by $\lambda_1(X)$ the smallest non-zero eigenvalue of $\Delta_X$, and often refer to it as
``the spectral gap". In the finite area situation, it may happen that there is no non-zero $L^2$-eigenvalue, in that case 
the bottom of the continuous spectrum is the relevant quantity and equals $1/4$.

\subsubsection{Link with geodesic counting.}

The spectral gap of hyperbolic surfaces is a ubiquitous quantity that appears in several natural
problems. One striking example of this is the question of counting closed geodesics up to a certain
length, historically pioneered by Heinz Huber \cite{huber1960} after the seminal work of Selberg on
the trace formula \cite{SelbergTrace}. Let us denote by $\mathcal{P}_X$ be the set of oriented primitive
closed geodesics on $X$, and if $\gamma \in \mathcal{P}_X$, $\ell(\gamma)$ will denote its
length. Then it is known that as $T\rightarrow +\infty$,
$$N(T):=\#\{ \gamma \in \mathcal{P}_X\ :\ \ell(\gamma)\leq T\} \sim \frac{e^T}{T},$$
a result called the prime number theorem due to its similarity with number theory. The leading order
$e^T/T$ is related to the trivial eigenvalue $0$. The spectral gap appears when studying the error
term in this approximation, and we have the more detailed asymptotic:
\begin{equation}
  \label{eq:prime}
  N(T)= \mathrm{li}(e^T)+\sum_{j} \mathrm{li}(e^{\alpha_j T})+O( e^{(3/4)T}),
\end{equation}
where $\mathrm{li}(x)=\int_2^x \frac{dt}{\log(t)}$ and the exponents
$1/2\leq\alpha_j\leq\ldots\leq\alpha_1<1$ are in correspondence with eigenvalues
$0<\lambda_1\leq \ldots \leq \lambda_j$ of $\Delta_X$ inside $[0,1/4]$ via the identity
$\alpha_j(1-\alpha_j)=\lambda_j$.  See for example \cite{Iwaniec} and references therein for more
detailed asymptotic formulas in the compact and finite volume case. Besides geodesic counting, the
spectral gap appears also naturally in the error term of the Hyperbolic lattice counting point
problem \cite[Chapter 12]{Iwaniec}, the rate of mixing of the geodesic flow on the unit tangent
bundle of $X$, see Ratner \cite{Ratner}, as well as the Brownian motion by work of Golubev--Kamber
\cite{golubev2019}.

\subsubsection{Link with connectivity.}
\label{s:connectivity}

Another very interesting feature of the spectral gap is that it is a quantitative measurement for
connectivity. The heuristic is that the trivial eigenvalue $0$ has multiplicity when the surface is
disconnected; hence, a ``poorly connected'' surface $X$ will have a small $\lambda_1(X)$. This can
be made precise by comparing $\lambda_1(X)$ with the \emph{Cheeger constant}
\begin{equation*}
  h(X) = \inf_A \frac{\ell(\partial A)}{\mathrm{Area}(A)}
\end{equation*}
where the infimum is taken over all smooth subsurfaces $A \subset X$ with $\mathrm{Area}(A)
\leq \frac12 \mathrm{Area}(X)$. We have the inequalities
\begin{equation*}
  \frac{h(X)^2}{4} \leq \lambda_1(X) \leq 2h(X)(h(X)+5)
\end{equation*}
the left part due to Jeff Cheeger \cite{cheeger1970} and right to Peter Buser \cite{buser1982}. In
other words, the spectral gap $\lambda_1(X)$ is close to zero if and only if the Cheeger constant
$h(X)$ is close to zero.  Another measurement of connectivity is the \emph{diameter}
$\mathrm{diam}(X)$, which measures the maximal distance between two points; Michael Magee explains
in \cite{magee2020b} how the diameter and $\lambda_1$ interplay by a interesting argument related to
the speed of mixing.

\subsubsection{Large-volume limit.}

A natural question in spectral geometry is the behaviour of $\lambda_1(X)$ in the regime of large
volume. For simply connected domains $\Omega \subset \H$, Osserman proved in \cite{Osserman} that
when looking at the first Dirichlet eigenvalue of $\Omega$, we have
$$\lambda_1(\Omega)\geq \frac{1}{4(\tanh(R(\Omega)))^2} > \frac 14$$
where $R(\Omega)$ is the in-radius i.e. the radius of largest inscribed geodesic discs in
$\Omega$. Therefore for hyperbolic domains, $\lambda_1(\Omega)$ is allowed to decrease to $1/4$ only if
$R(\Omega)$ goes to infinity. Here the number~$1/4$ appears as the bottom of the spectrum of
the universal cover $\Delta_{\mathbb{H}^2}$.

The situation for compact hyperbolic surfaces is completely different: Huber proved in
\cite{huber1974} that
\begin{equation}
  \label{eq:huber_14}
  \limsup_{k \rightarrow \infty} \lambda_1(X_k) \leq \frac 14
\end{equation}
for any sequence $(X_k)_k$ of compact hyperbolic surfaces of volume going to infinity. Note that, in
this regime, by the Gauss--Bonnet formula, the volume is equal to $4 \pi (g-1)$, where $g$ is the
genus of the surface; hence, the volume and genus go to infinity hand in hand.

Non-zero eigenvalues under $1/4$ for finite area surfaces are somewhow
``special'' because they live outside the bulk spectrum of the universal cover; for instance, they
are related to subdominant terms in \eqref{eq:prime}. For this reason, they are often referred to as
\emph{small eigenvalues}.

Even though McKean erronously affirmed that such eigenvalues cannot exist \cite{mckean1972}, it is
fairly easy to produce surfaces~$X$ with large volume and arbitrarily small $\lambda_1(X)$ by
looking at large cyclic covers over a fixed base surface: see Randol \cite{Randol} for a proof based
on abelian twists and Selberg's trace formula.  For a more geometric proof based on Cheeger--Buser's
inequality yiedling examples in any genus, see in Bergeron \cite[Chapter 3]{BergeronBook}. Actually,
this geometric argument can be used to produce hyperbolic surfaces of signature $(g,n)$ with up to
$2g-2+n$ eigenvalues below $\epsilon$ (where $\epsilon > 0$ can be taken arbitrarily small), by
pinching a maximal family of simple disjoint closed curves. We know this to be the maximal number of
small eigenvalues a finite-area surface can have by work of Otal and Rosas \cite{otal2009}.

By analogy with cubic graphs, Buser mistakenly conjectured in \cite{buser1978} that $\lambda_1(X)$
should always tend to zero as $\mathrm{Vol}(X)$ tends to infinity, so that the $1/4$ in
\eqref{eq:huber_14} can be replaced by a $0$. On the other hand, Selberg
famously conjectured the opposite behaviour for congruence surfaces.

\begin{conjecture} (Selberg \cite{selberg1965})
Let $\Gamma(N)$ denote the discrete subgroup of $\mathrm{PSL}_2(\Z)$ defined by
$$\Gamma(N):=\{ g \in  \mathrm{PSL}_2(\Z)\ :\ g \equiv \pm Id \ \mathrm{mod}\ N\},$$
and set $X(N):=\Gamma(N)\backslash \H$. Then for all $N\geq 1$ one has $\lambda_1(X(N))\geq 1/4$.
\end{conjecture}
Notice that $X(N)\rightarrow X(1)$ is a Galois cover and we have $\mathrm{Vol}(X(N))=\vert \mathrm{PSL}_2(\Z_N)\vert \mathrm{Vol}(X(1))$, thus Selberg's eigenvalue conjecture is really about
understanding $\lambda_1$ in the large volume regime. In the same paper \cite{selberg1965}, Selberg proved that for all $N\geq 1$ we have $\lambda_1(X(N))\geq 3/16$. He also mentioned,
by looking at abelian covers, that his conjecture was wrong when considering arbitrary finite index
subgroups. Selberg's eigenvalue conjecture has major implications in classical number theory and
part of the bigger picture of ``Ramanujan's conjectures''. We refer the reader to \cite{Sarnak1} for
more insights on the importance and difficulty of this conjecture. Selberg's lower bound has since
been improved, to our knowledge the best bound so far is due to Kim and Sarnak \cite{KimSarnak} and
is $\lambda_1(X(N))\geq 975/4096 \approx 0.2380$.

In \cite{buser1988}, Buser, Burger and Dodziuk, based on the above Selberg bound, showed how to
construct a sequence $(X_k)_k$ of compact hyperbolic surfaces with genus (and thus volume) going to
infinity while $\limsup_k \lambda_1(X_k)\geq 3/16$. In the very same paper, they conjectured that
one should be able to do the same with $1/4$ in place of $3/16$. More generally, compact arithmetic congruence surfaces do satisfy $\lambda_1\geq 3/16$, see the work of
M.F. Vigneras \cite{Vigne1}.

\subsubsection{Benjamini--Schramm convergence.} A way to formalize the fact that a sequence $(X_k)$ of compact hyperbolic surfaces with genus going to infinity looks increasingly like the hyperbolic plane $\H$ is through the notion of Benjamini--Schramm convergence, which comes from the world of graphs \cite{BS1}. More preciselly, let $\mathrm{Inj}_z(X)$ denote the injectivity radius at $z$ on $X$. Assume that
$\lim_{k\rightarrow \infty} \mathrm{Vol}(X_k)=+\infty$. We say that $(X_k)_k$ converges in the sense of Benjamini--Schramm to $\H$ if we have for all $L>0$,
$$\lim_{k\rightarrow +\infty} \frac{\mathrm{Vol}(\{z\in X_k\ :\ \mathrm{Inj}_z(X_k)<L\})}{\mathrm{Vol}(X_k)}=0. $$
In other words, when picking a point $z$ randomly on $X_k$ according to hyperbolic volume, the
probability that the ball centered at $z$ of radius $L$ is isometric to a ball of same radius in
$\H$ goes to $1$  in the the large volume limit. The systole $\ell(X)$ of a surface $X$ is the shortest length of closed geodesics. A sequence of surfaces
$(X_k)$ with genus going to infinity is called "uniformly discrete" if there exists a positive uniform lower bound on $\ell(X_k)$.

Fix $a<b$ in $\R$, and let $N_k(a,b)$ denote the number of eigenvalues of the Laplacian on $X_k$
inside the interval $[a,b]$, counted with multiplicity.  It is shown by Le Masson and Sahlsten in
\cite{LS1}, see also \cite{monkapde} for a more quantitative version, that under the assumption of
Benjamini--Schramm convergence of $(X_k)_k$ to $\H$ plus uniform discreteness, one always has the limit
$$\lim_{k\rightarrow +\infty} \frac{N_k(a,b)}{\mathrm{Vol(X_k)}}=\frac{1}{4\pi}\int_{1/4}^{+\infty} \mathbf{1}_{[a,b]}(\lambda)\tanh\left(\pi\sqrt{\lambda-\frac{1}{4}}\right)d\lambda. $$
Therefore, under those two assumptions, the spectral density (normalized by volume) of the
Laplacian on $X_k$ converges to the spectral measure of the Laplacian on $\H$, supported on
$[1/4,+\infty)$. We point out that Benjamini--Schramm convergence (without uniform discreteness) does not implies spectral convergence, see
\cite{GK}. The notions of Benjamini--Schramm convergence and spectral convergence have a far reaching generalization for sequences of higher dimensional manifolds: we mention here the work of the "seven samurais" \cite{seven1} and references herein.

It seems plausible that the typical spectral gap of a large hyperbolic surface should be close to
$1/4$. However Benjamini--Schramm convergence alone is {\it not enough} to get a uniform spectral
gap. Indeed, it is possible to build sequences $(X_k)_{k\in \N}$ of compact surfaces
with genus $g_k$ going to infinity as $k\rightarrow +\infty$ such that $\inf_{z\in X_k} \mathrm{Inj}_z(X_k)\rightarrow \infty$ while $\lambda_1(X_k)=O(\frac{\log(g_k)}{g_k})$. This construction follows from Buser--Sarnak \cite{BS2}, where they show that compact arithmetic surfaces $S_g$ with genus $g$ have systoles of size $\log(g)$ and hence their injectivity radius goes to infinity with $g$. On the other hand, they prove in \cite{BS2} that it is always possible to find a simple closed geodesic $\gamma$ on $S_g$ which is {\it non trivial in homology} and such that $\ell(\gamma)\leq 2\log(4g-2)$. By cutting and pasting along this geodesic, one can then build a connected $2$-cover $\widetilde{S_g}$ of $S_g$ with Cheeger constant of size at most $O(\log(g)/g)$ and therefore $\lambda_1(\widetilde{S_g})$ becomes arbitrarily small as $g\rightarrow \infty$.
These facts highlight the intuition that to be able to prove a near optimal uniform spectral gap, probabilistic tools might be helpful. The existing huge literature on random graphs provides some insight on the matter.

\subsection{Random regular graphs are optimal expanders.}

Let us make a short scenic detour into the realm of random regular graphs where many analogous
questions have been settled.

Let $\mathcal{G}=\mathcal{G}(V,E)$ be a (non-oriented) $k$-regular graph with finite set of edges
$E$ and vertices $V$. In what follow we will set $n=\# V$ and we will be interested only in the regime
where $k$ remains fixed and $n$ goes to infinity. Let $A$ denote the adjacency matrix of $\mathcal{G}$, i.e.  the $n\times n$ symmetric matrix defined by 
$A(i,j):=\#\{\mathrm{edges\ connecting\ }i\ \mathrm{and}\ j\}$. Because $\mathcal{G}$ is $k$-regular, the spectrum of $A$ consists of $n$ real eigenvalues included in $[-k,k]$ with $k$ being simple if and only if $\mathcal{G}$ is connected. We define the {\it spectral gap of connected graphs} by 
$$\lambda(\mathcal G)=\max\{|\lambda| : \lambda\ \text{ is an eigenvalue of }A\
\mathrm{and}\ \lambda\neq \pm k\}.$$ This spectral gap $\lambda(\mathcal{G})$ is a natural quantity
related to the speed of convergence to equilibrium of the Markov chain $\frac{1}{k }A$.

There is a lower bound to $\lambda(\mathcal{G})$ for $k$-regular graphs in the large $n$ regime:
more precisely, if $\mathcal{G}_{n,k}$ is a family of connected $k$ regular graphs, then
$$\liminf_{n\rightarrow +\infty} \lambda(\mathcal{G}_{n,k})\geq 2\sqrt{k-1}.$$
This is the Alon--Boppana inequality, see in \cite{LPS} for a proof. The quantity $2\sqrt{k-1}$
appears here as the spectral radius of the uniform random walk on the infinite $k$-regular tree,
which was first observed by Kesten \cite{Kesten}. Therefore $2\sqrt{k-1}$ is the perfect analog of
$1/4$ in the hyperbolic surface setting.

A $k$-regular graph is called {\it Ramanujan} if $\lambda(\mathcal{G})\leq 2\sqrt{k-1}$. Explicit
families of $k$-regular expander graphs for specific values of $k$ where constructed in \cite{LPS}
by subtle number theoretic methods. It is somehow surprising how elusive these Ramanujan graphs are:
indeed, Alon famously conjectured in 1986 that ``most'' large $k$-regular graphs are
almost-Ramanujan \cite{alon1986}. To address this conjecture, Broder and Shamir studied the
concentration of the second largest eigenvalue $\lambda_2(\mathcal{G})$ of random $k$-regular
graphs, and proved, for different models, a probabilistic bound of the form
$\lambda_2(\mathcal{G}) \leq 3 k^{3/4}$ \cite{broder1987}. The first proof of the Alon conjecture
was obtained by Friedman in 2003, who showed in \cite{friedman2003} that, for all $\epsilon >0$,
\begin{equation*}
  \lim_{n \rightarrow \infty} \mathbb{P}_{k,n}(\lambda(\mathcal{G}) \leq 2 \sqrt{k-1}+\epsilon)
  = 1
\end{equation*}
for various usual probability distributions $\mathbb{P}_{k,n}$ on the set of $k$-regular graphs with
$n$ vertices. The proof is a very technical application of the trace method, with a subtle argument
allowing to discard problematic patterns (called ``tangles'').  Many years later, Bordenave has
provided a new much simpler proof of this statement using the non-backtracking matrix
\cite{bordenave2020}.

A similar line of intuitions and results can be found in the study of random lifts. If
$\mathcal{G} = \mathcal{G}(V,E)$ is a starting $k$-regular graph, we can construct a random lift
$\mathcal{G}_p = (V_p,E_p)$ of $\mathcal{G}$ of degree $p$ by letting
$V_p := V \times \{1, \ldots, p\}$ and, for each edge $(i,j) \in E$, picking a random bijection from
$\{i\} \times \{1, \ldots, p\}$ to $\{j\} \times \{1, \ldots, p\}$ to constitute the edges
$E_p$. Then one can recover the initial graph $\mathcal{G}$ by forgetting the second component,
i.e. the map $\pi : \mathcal{G}_p \rightarrow \mathcal{G}$ is a covering with fiber of size $p$.

In this setting, the spectrum of the adjacency matrix $A^{(p)}$ of $\mathcal{G}_p$ contains a copy
of the spectrum of the adjacency matrix $A$ of $\mathcal{G}$, because any eigenfunction on
$\mathcal{G}$ can be lifted to be an eigenfunction on $\mathcal{G}_p$ constant on each fiber. The
interesting question now becomes the \emph{new eigenvalues}, and the corresponding spectral gap
\begin{equation*}
  \lambda_{\mathrm{new}}(\mathcal{G}_p)
  = \max \{ |\lambda| : \ \lambda \text{ is an eigenvalue of } A^{(p)} \text{ and not of } A \}.
\end{equation*}
Friedman conjectured in \cite{friedman2003a} that, for a typical lift of large degree $p$,
$\lambda_{\mathrm{new}}(\mathcal{G}_p) \leq 2 \sqrt{k-1} + \epsilon$, i.e. that Alon's conjecture is
also true of random lifts. This was proven by Friedman and Kohler in \cite{friedman2019} and
improved shortly afterwards by Bordenave in his shorter proof
\cite{bordenave2020}. In the same line of probabilistic ideas, Bordenave and Collins proved strong convergence for sums of random permutations matrices \cite{bordenave2019},
a key tool that will prove also to be critical in the study of random covers of hyperbolic surfaces.

Recent years have been marked by two breakthroughs in this topic, showing that this is still a very
active question. On the one hand, Huang, McKenzie and Yau have proved in \cite{huang2025} that, for
large enough $n$, the probability for a random $k$-regular graph with $n$ vertices to be Ramanujan
is approximately $69 \%$. This builds on significant new developments in the study of Green
functions. On the other hand, Chen, Garza-Vargas, Tropp and van Handel provided a
surprising new approach, distinct from trace methods and Green functions, allowing to establish
strong convergence and hence dramatically simplify the proof of Friedman's theorem
\cite{chen2025}. We will see how trace methods work and discuss the notion of strong convergence in
more detail when we come back to the study of hyperbolic surfaces.

\subsection{Models of random hyperbolic surfaces and their spectral gap.}

It is now time to come back to the initial question, the spectral gap of hyperbolic surfaces of
finite area. Regular graphs and hyperbolic surfaces often behave in similar fashion, as homogenous
spaces with exponential balls (the volume of a ball of radius $R$ in the infinite $k$-regular tree
and the hyperbolic plane $\mathbb{H}$ are $k(k-1)^{R-1}$ and $2 \pi (\cosh R-1)$ respectively). The
intuition coming from graphs and Selberg's eigenvalue conjecture leads naturally to the following
{\it meta-conjecture}.

\begin{conjecture}
  Sequences of random hyperbolic surfaces $(X_k)_k$, with
  $\lim_{k\rightarrow \infty}\mathrm{Vol}(X_k)= +\infty$, have a near optimal spectral gap
  with high probability, i.e for all $\epsilon>0$,
 $$\lim_{k\rightarrow +\infty} \mathbb{P}\left ( \lambda_1(X_k)\geq \frac{1}{4}-\epsilon \right)=1.$$
\end{conjecture}

Note that, by ``with high probability'', we mean with probability going to $1$ as the asymptotic
parameter goes to infinity.  Now this raises a big question: how does one study random hyperbolic
surfaces? Whilst the set of $k$-regular graphs with $n$ vertices is finite and can be equipped
naturally with a variety of interesting combinatorial probability measures, the space of all
hyperbolic metrics is much richer.  In this survey, we will present many different probabilistic
settings, each with their own strengths and specific interest. Some are continuous (most notably the
Weil--Petersson measure on moduli spaces) whilst others are discrete (the Brooks--Makeover model and
random covers, as defined below).  In the last few years, we have reached a precise understanding of
the spectral gap of typical hyperbolic surfaces, either sampled with respect to the Weil--Petersson
volume, or as uniform random covers (which shall be reminiscent of the case of graphs). The aim of
these proceedings is to present these advances and how they relate with one another and the parallel
world of graphs.

\section{First families of expanders: the Cheeger method.}

Let us present the first wave of results regarding the spectral gap of large hyperbolic surfaces
using probabilistic methods. The unifying factor here is that the lower bound on the spectral gap
$\lambda_1$ is obtained by obtaining a probabilistic lower bound on the Cheeger constant $h$
introduced in \ref{s:connectivity}, together with the Cheeger inequality $\lambda_1 \geq h^2/4$.

\subsection{The Brooks--Makover model.}
\label{sec:brooks-makover-model}

The first probabilistic model allowing to study random hyperbolic surfaces was a combinatorial
construction suggested by Brooks and Makover \cite{brooks2004}. The idea is to gather a collection
of $2n$ oriented ideal hyperbolic triangles (i.e. triangles with all angles equal to $0$), and glue
them (with no shear) following a random oriented $3$-regular graph. The resulting surface is a
random hyperbolic surface of area $2\pi n$ of genus strongly concentrated around $n/2$, with cusps
corresponding to specific cycles in the random graph. One can compactify this surface and hence
obtain a random compact hyperbolic surface of typically large genus. The surfaces obtained by this
construction are exactly Bely\u{\i} surfaces, i.e. the Riemann surfaces which can be constructed
over some number field \cite{belyi1979}; in particular, they are dense in the set of all hyperbolic
surfaces.  Comparing the geometry of the hyperbolic surface to its embedded $3$-regular graph, and
using the Cheeger inequality, Brooks and Makover proved the following.

\begin{theorem}[{\cite{brooks2004}}]
  There exists $C>0$ such that $\lambda_1(X) \geq C$ with high probability  for the
  Brooks--Makover model.
\end{theorem}

The value of the universal constant $C$ is not explicit (although it must be $\leq 1/4$ by Huber's
bound \ref{eq:huber_14}). Interestingly, Shen and Wu proved that, for any $\epsilon >0$,
$h(X) \leq \frac{3}{2 \pi} + \epsilon$ with high probability for the Brooks--Makover model. The
important thing here is that the value $\frac{3}{2 \pi} \approx 0.4775$ is stricly smaller than the
Cheeger constant of the hyperbolic plane $h(\H)=1$, for which the two sides of the Cheeger
inequality match.  This means that there is no chance to reach the optimal spectral gap $1/4$ using
a probabilistic bound together with the Cheeger inequality. This observation has been extended by
Budzinski, Curien and Petri \cite{budzinski2022}, who proved that
  \begin{equation*}
  \limsup_{g \rightarrow \infty} \sup_{X \in \mathcal{M}_g} h(X) \leq \frac{2}{\pi} \approx 0.6366.
\end{equation*}
The argument relies on a random coloring of the surface, and is adapted from a similar result for
regular graphs proven by Bollob\'as~\cite{bollobas1988} and improved by Alon \cite{alon1997}.  It
follows that the spectral gap and Cheeger constant are expected to behave quite differently in the
large-volume limit, one approaching the universal cover whilst the other does not. However, because
the geometry of random hyperbolic surfaces is oftenwise easier to access than their spectrum,
bounding the Cheeger constant is often the first way to obtain a spectral gap result.

It is worth noting that another combinatorial construction using a random $3$-regular graph as a
skeleton, now gluing hyperbolic pairs of pants (with all boundary lengths fixed to a specific
parameter and no twists), was used by Budzinski--Curien--Petri \cite{budzinski2021} to find the
first examples of large hyperbolic surfaces with optimal diameter. This model was enriched to be
continuous by Mathien \cite{mathien2024}. However, by a simple min-max argument, the spectral gap of
these random surfaces is very small, as they contain pair of pants with a short boundary (see
e.g. \cite[Theorem 8.1.3]{buser1982}).

\subsection{The Weil--Petersson model.}

Beyond the Brooks--Makover combinatorial model, there is a very natural continuous probabilistic
model allowing to study random hyperbolic surfaces, which we will now introduce in details, as it is
one of the core focuses of these proceedings.

\subsubsection{The sample space.}

Let us fix a genus $g \geq 2$. The natural sample space we wish to consider is the \emph{moduli space}
\begin{equation*}
  \mathcal{M}_{g} := \{\text{hyperbolic surfaces of genus } g\}
  \diagup
  \{\text{isometries}\}.
\end{equation*}
Due to the famous Uniformisation Theorem by Klein--Poincaré--Koebe, we can alternatively define
\begin{equation*}
  \mathcal{M}_{g} := \{\text{Riemann surface of genus } g\} 
  \diagup
  \{\text{biholomorphism}\}
\end{equation*}
where a Riemann surface is a one-dimensional complex manifold. This is the definition that Riemann
used when introducing the moduli space $\mathcal{M}_g$ in 1857 \cite{riemann1857}. He subsequently
argued that it was a space of complex dimension $3g-3$ by interpreting the elements of
$\mathcal{M}_g$ as branched covering surfaces of the Riemann sphere and counting the number of free
parameters using the Riemann--Roch theorem. 

The general approach to study the moduli space $\mathcal{M}_g$ shifted over the course of the 20th
century to the study of its universal cover $\mathcal{T}_g$ called the \emph{Teichm\"uller
  space}. This is a very old field which has branches in many areas of mathematics; we refer to
\cite{imayoshi1992} for a more detailed presentation of the history below.

The difference between a point of the moduli space $\mathcal{M}_g$ and its universal cover
$\mathcal{T}_g$ is that an element $X$ of $\mathcal{T}_g$ comes with a \emph{marking}, which can,
for instance, be viewed as a fixed set of generators of the fundamental group $\pi_1(X)$ of
$X$. More precisely, for a genus $g \geq 2$, we define the \emph{surface group} of genus $g$ by the
presentation 
\begin{equation}
  \label{eq:pi_1}
  \Gamma_g
  := \langle (a_i, b_i)_{1 \leq i \leq g} \, | \, [a_1, b_1] \ldots [a_g,b_g] \rangle.
\end{equation}
A point in $\mathcal{T}_g$ is the data of an element $X$ of the moduli space $\mathcal{M}_g$,
together with a homomorphism between $\Gamma_g$ and $\pi_1(X)$, defined up to an inner automorphism.

To recover the moduli space from its universal cover, we need to ``forget the marking''. This is
done by taking the quotient of $\mathcal{T}_g$ by the action of the mapping class group
\begin{equation*}
  \mathrm{MCG}_g := \{\text{positive diffeomorphisms of a surface of genus } g\}
  \diagup
  \{\text{isotopy}\}.
\end{equation*}
This is a discrete group acting properly discontinuously on $\mathcal{T}_g$. We shall see very soon
that $\mathcal{T}_g$ is a manifold. Because the action of $\mathrm{MCG}_g$ is not free and some
elements have fixed points, the moduli space $\mathcal{M}_g$ is not a manifold but rather an
orbifold, of the same dimension as $\mathcal{T}_g$.

\subsubsection{Parametrizations.}
\label{sec:param-space-hyperb}

It is useful when studying a space to have a description in terms of coordinates.  Several
parametrizations of $\mathcal{T}_g$ in terms of nice, geometric coordinates can be obtained, placing
oneself on the hyperbolic side of the Uniformization Theorem.  The first one consists in identifying
$\mathcal{T}_g$ with the space of discrete faithful representations of the surface group~$\Gamma_g$
in $\mathrm{PSL}_2(\R)$ (the group of positive isometries of $\mathbb{H}^2$). This approach yields a
system of $6g-6$ Fricke coordinates, which first appear in a 1897 book by Fricke and
Klein~\cite{fricke1897}.

\begin{figure}[htbp]
  \centering
  \includegraphics[scale=0.4]{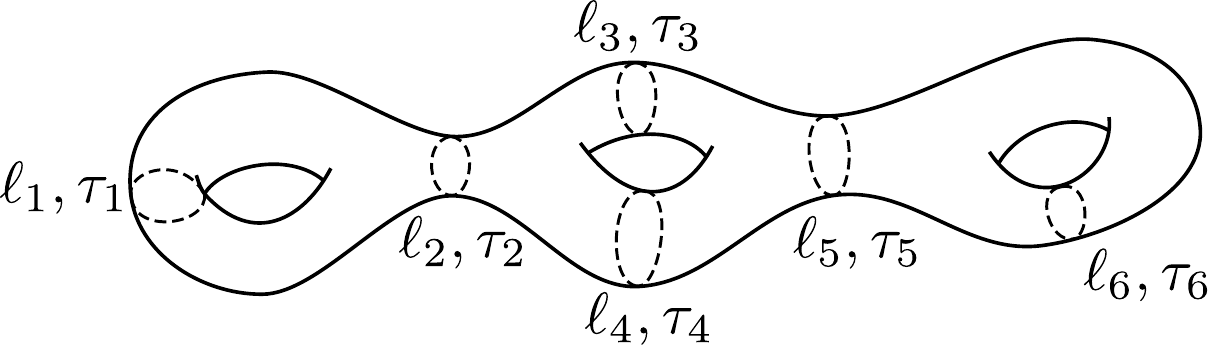}
  \caption{Fenchel--Nielsen coordinates in genus $g=3$.}
  \label{fig:FN}
\end{figure}

A very important set of coordinates for us is the one introduced by Fenchel and Nielsen in an
unpublished manuscript in 1948.  These coordinates are constructed by considering a decomposition of
our hyperbolic surface in hyperbolic pairs of pants, i.e. surfaces of genus $0$ with three geodesic
boundary components, as represented in \cref{fig:FN}. By additivity of the Euler characteristic,
such a decomposition always contains $2g-2$ pairs of pants, delimited by $3g-3$ curves. The (marked)
geometry then is entirely determined by the data of the lengths $(\ell_i)_{1 \leq i \leq 3g-3}$ of
the boundary components of the pairs of pants, together with the \emph{twist parameters}
$(\tau_i)_{1 \leq i \leq 3g-3}$ describing the gluings of the various pairs of pants.
In the moduli space, the twist parameter $\tau_i$ can be viewed as an element of
$\R \diagup \ell_i \R$, whilst in the universal cover $\mathcal{T}_g$ it is an element of $\R$.

\subsubsection{The Weil--Petersson symplectic form.}
\label{sec:weil-petersson-form}

The next step in our exploration of the moduli space $\mathcal{M}_g$ is to equip it with a
probabilistic measure.  In order to do so, we shall continue the back-and-forth and we go back to
the complex-analysis side of the Uniformisation Theorem.

The story starts by one question: how do we deform a complex structure? Or, in another words, is
there any nice description of the tangent spaces of $\mathcal{T}_g$ and $\mathcal{M}_g$? It was
proven, heuristically by Teichm\"uller in the 1940s and formally by Ahlfors and Bers in the 1950s,
that the cotangent space at a point $X \in \mathcal{T}_g$ can be identified with the space $A_2(X)$
of \emph{quadratic differentials} on the Riemann surface $X$.\footnote{Another complementary
  viewpoint is to use the more general Kodaira--Spencer theory on the deformation of complex
  structures (as preferred, notably, by Weil when introducing the Weil--Petersson form in
  \cite{weil1958a}).} The important point here is that the space $A_2(X)$ has a natural inner
product called the \emph{Petersson inner product} (which is, basically, the $L^2$ inner product on
$X$, but for modular forms). In 1958, Weil suggested to study the space $\mathcal{T}_g$ through the
lense of this metric, and raised several open questions \cite{weil1958a}. Ahlfors soon followed up
and answered several of Weil's questions, notably proving that this metric was K\"ahler
\cite{ahlfors1961a}. In particular, we obtain a symplectic form called the
\emph{Weil--Petersson form} $\omega_g^{\mathrm{WP}}$ on $\mathcal{T}_g$, which is naturally
$\mathrm{MCG}_g$-invariant and hence descends onto the moduli space~$\mathcal{M}_g$.

It is worth noting that, at this point in time (i.e. in the 1960s), we have everything required to
define a probabilistic model on the moduli space~$\mathcal{M}_g$. Indeed, any symplectic form
induces a volume form, given by
$$\dvwp := \frac{1}{(3g-3)!} (\omega_g^{\mathrm{WP}})^{\wedge(3g-3)},$$ where $3g-3$ is half the
dimension of the sympleptic space. Since the total volume of the moduli space $\mathcal{M}_g$ is
finite, it can be renormalized to obtain a probability measure $\pwp$ on $\mathcal{M}_g$. However,
beyond this abstract definition, no tools are yet available to allow the study of surfaces sampled
following $\pwp$ (and we will see there are multiple challenges ahead!). It is unclear whether this
question was left out during 50 years due to lack of interest in probabilistic methods in favor of
dynamical or geometric questions, or these technical challenges; probably a bit of both.

\subsubsection{Wolpert and Mirzakhani's magic formulae.}
\label{sec:wolp-mirz-magic}

In the early 2010s, Mirzakhani and Guth--Parlier--Young independently have the idea to use the
Weil--Petersson volume form to study random hyperbolic surfaces, notably in the large-volume limit
\cite{guth2011,mirzakhani2013}. A cornerstone of both articles is a formula due to Wolpert, often
referred to as \emph{Wolpert's magic formula} \cite{wolpert1981}, which states, very simply, that
for \emph{any} set of Fenchel--Nielsen coordinates $(\ell_i, \tau_i)_{1 \leq i \leq 3g-3}$,
\begin{equation*}
   \omega_g^{\mathrm{WP}} = \sum_{i=1}^{3g-3} \d \ell_i \wedge \d \tau_i.
\end{equation*}
In other words, Fenchel--Nielsen parameters are sympleptic coordinates for the Weil--Ptersson form
and, in particular, the Weil--Petersson volume $\dvwp$ is the Lebesgue measure in these coordinates.
The miracle in this formula is that it allows to connect the two sides of the Uniformization
Theorem: whilst the definition of the Weil--Petersson form relies deeply on the complex structure,
Fenchel--Nielsen coordinates are purely hyperbolic quantities. Another striking feature is that the
formula holds for any set of Fenchel--Nielsen coordinates, which is non-trivial, especially given
there exists distinct decompositions of a surface of genus $g$ which yield distinct sets of
coordinates, and are hard to relate to one another. A modern reformulation of this result can be
stated as follows.

\begin{theorem}[Wolpert, \cite{wolpert1981}]
  \label{thm:wolpert}
  The symplepctic form $\sum_{i=1}^{3g-3} \d \ell_i \wedge \d \tau_i$ on the Teichm\"uller space
  $\mathcal{T}_g$ is independent on the choice of a Fenchel--Nielsen set of coordinates and
  $\mathrm{MCG}_g$-invariant, and hence descends to a uniquely defined sympleptic form on the moduli
  space $\mathcal{M}_g$.
\end{theorem}

In other words, our space of interest, the moduli space $\mathcal{M}_g$, is ``simply'' the quotient
of $(\R_{>0} \times \R)^{3g-3}$ by a discrete group, and equipped with the renormalized Lebesgue
measure. This seems like a pretty simple picture, and indeed this is all that Guth, Parlier and Young
need to prove probabilistic bounds on the lengths of any pair of pants decomposition for typical
large random hyperbolic surfaces sampled with $\pwp$, one of the first geometric information
obtained on these objects \cite{guth2011}.

There is however a great challenge one needs to tackle in order to move forward and provide a
precise description of random surfaces sampled using $\pwp$. In order to integrate on a quotient
space (think for instance of the torus $\R^d \diagup \Z^d$), it is essential to understand the
quotient well. The integration on the quotient space is oftenwise performed by integrating over a
fundamental domain ($[0,1]^d$ in the example). Unfortunately, the action of the mapping class group
on $\mathcal{T}_g$ is very delicate to grasp. We do not know any explicit fundamental domain for
this action, which severely limits our understanding of the space. One can construct a rough
fundamental domain (i.e. a subset of $\mathcal{T} _g$ containing at least one representant of each
$\mathrm{MCG}_g$-orbit, but with some overlap) using a result by Bers~\cite{bers1985}, which states
there exists a constant $B_g$ such that any hyperbolic surfaces admits a decomposition in pairs of
pants will all lengths $\leq B_g$ (see \cite[Proposition 3.5]{monk2021} for an upper bound on the
total volume of $\mathcal{M}_g$ using this approach). However this is a  rough description and
not enough to develop a theory of integration over~$\mathcal{M}_g$.

This issue was completely circumvented by using a radically different approach, relying on two
magical formulae due to Mirzakhani.  On the one hand, \emph{Mirzakhani's integration formula}
\cite{mirzakhani2007} allows to compute the integral over the moduli space of certain specific types
of functions, called \emph{geometric functions}, i.e. functions of the form
$\sum_\gamma F(\ell(\gamma))$, where $F : \R_{\geq 0} \rightarrow \R$ is a positive measurable
function, and $\gamma$ runs over a $\mathrm{MCG}_g$-invariant set of simple closed geodesics. The
proof of this formula consists in identifying the level-sets for the functions $\ell(\gamma)$ with
moduli spaces (of surfaces with a geodesic boundary) through a symplectomorphism.  On the other
hand, \emph{Mirzakhani's partition of unity} \cite{mirzakhani2004} is a generalization of the famous
identity due to McShane \cite{mcshane1991} which states that, for \emph{any} hyperbolic metric on
the once-punctured torus,
\begin{equation}
  \label{eq:mcshane}
   \sum_{\gamma} \frac{1}{1+e^{\ell(\gamma)}} = \frac 12
\end{equation}
where the sum runs over all unoriented simple closed geodesics $\gamma$ on the once-punctured torus.
The brilliant idea here is to view this as a rewriting of the constant function equal to $1/2$ (the
r.h.s.)  as a function of the lengths of simple closed geodesics (the l.h.s.). The left-hand side of
\cref{eq:mcshane} is a combination of geometric functions that can be integrated using Mirzakhani's
integration formula; this leads to an expression for the total volume of the moduli space of the
once-punctured torus. More generally, Mirzakhani's partition of unity allows to express the constant
function as a combination of geometric functions, and hence find a recursive formula for
Weil--Petersson volumes of moduli space.

\subsubsection{Probabilistic results in this model.}

Mirzakhani's integration formula, together with the recursive formula derived from the partition of
unity, provide a very exciting framework to study the geometry of random hyperbolic surfaces sampled
with $\pwp$.  In \cite{mirzakhani2013}, Mirzakhani presents a full proof of concept as well as
multiple new geometric results, including the first bound on the spectral gap of typical hyperbolic
surfaces of high genus
\begin{equation*}
  \lim_{g \rightarrow \infty}
  \mathbb{P}_g^{\mathrm{WP}}\left(\lambda_1 \geq \frac 14 \left(\frac{\log 2}{\log 2+\pi}\right)^2\right)=1.
\end{equation*}
As in the case of the Brooks--Makover model, the value
$\left(\log 2/(\log 2+\pi)\right)^2/4 \approx 0.002$ comes from the Cheeger inequality and is a
priori not optimal. This is because the integration formula only allows to integrate geometric
functions, and hence it is much more manageable to study the Cheeger constant than the spectral
gap. Obtaining a uniform spectral gap for such a wide class of large-volume hyperbolic surfaces is a
great achievement and opens an exciting new research direction.

\section{Trace methods: to intermediate and optimal spectral gaps.}
\label{sec:trace-meth-refin}

Let us now discuss a second wave of spectral-gap results, improving on rough Cheeger-type
inequalities to obtain better (and even optimal) spectral gap results. We shall present a few ideas
behind trace methods in general, and then present results in the setting of random covers and
Weil--Petersson surfaces.

\subsection{A word on trace methods.}

The common technique shared by the results presented in this section is the \emph{trace method}. For
graphs, the trace method consists in studying the spectral gap through the asymptotic behaviour of
the traces of polynomials of the adjacency matrix. In the context of hyperbolic surfaces, the
situation is very much the same, now relying on the \emph{Selberg trace formula}
\cite{SelbergTrace}. We refer to \cite{anantharaman2024} for a detailed explanation of the following
method and here restrict the presentation to a few important points.

For a compact hyperbolic surface $X$ with eigenvalues $\lambda_j(X) = \frac14 + r_j(X)^2$, the
Selberg trace formula applied to an appropriate test function $h$ (e.g. a smooth even function
$\R \rightarrow \R$ with compact support) reads:
\begin{equation}
  \label{eq:stf}
  \sum_{j=0}^\infty \hat{h}(r_j(X)) = (g-1) \int_\R \hat{h}(r) r \tanh(\pi r) \d r
  + \sum_{\gamma \in \mathcal{P}_X} \sum_{k=1}^\infty
  \frac{h(k \ell(\gamma))}{2 \sinh (k \ell(\gamma)/2)}
\end{equation}
where we recall that $\mathcal{P}_X$ denotes the set of primitive oriented closed geodesics on $X$,
and $\hat h$ is the Fourier transform of $h$. This expression can be written differently in the
special case where $X$ is a cover, and generalised to finite-area surfaces by incorporating
additional terms.

Using the Selberg trace formula is a game of picking a good test function (or family of test
functions), which if well-curated can help us describe several aspects of the interactions between
the spectrum $(\lambda_j(X))_{j \geq 0}$ and the lengths of the closed geodesics
$\gamma \in \mathcal{P}_X$.  A simple way to make the small eigenvalues $\lambda_j (X) < \frac 14$
pop up in the formula above is the choice of a dilated test functions $h_L(r) := h(r/L)$ for a
function $h$ supported on $[-1,1]$. The parameter $L$ is a length cutoff: because $h_L$ is supported
on $[-L,L]$, the r.h.s. of \eqref{eq:stf} is restricted to closed geodesics of length $\leq L$. On
the other side of the formula, the biggest term in the l.h.s. is the one associated to the trivial
eigenvalue $\lambda_0=0$, which is equal to
\begin{equation}
  \label{eq:zero_eig}
  \hat{h}_L(i/2)
  = 2 \int_0^L h_L(r) \cosh \Big( \frac{\ell}{2} \Big) \d \ell
  = 2 L \int_0^1 h(r) \cosh \Big(\frac{\ell L}{2}\Big) \d \ell
  \gg e^{(1-\epsilon)L/2} \quad \text{for } L \gg 1.
\end{equation}
This exponential growth is matched on the other side of the equation by the fact that the number of
closed geodesics of length $\leq L$ behaves like $e^L$ as $L \gg 1$.

Looking at other terms, for any $j$, the Fourier transform $\hat{h}_L(r_j(X))$ will behave
exponentially in $L$ if and only if $r_j(X) \notin \R$, and be oscillatory otherwise. In other
words, the small eigenvalues $0 < \lambda_j(X) < 1/4$ stand out as the exponential terms in the
Selberg trace formula, their exponent being bigger the closer they are to the trivial eigenvalue $0$
(this should be reminiscent to the prime number theorem with error terms \eqref{eq:prime}).  The
general idea behind the trace method for graphs or surfaces alike consists in \emph{estimating the
  average of the trace formula}, and proving that, \emph{on average}, there are no exponential
constributions outside that of the trivial eigenvalue (or that we can at least bound their
exponential rates). This then allows to conclude that, \emph{with high probability}, there are no
small eigenvalues (or at least no non-trivial eigenvalues too close to~$0$).  The quality of the
spectral gap obtained depends on the precision of our estimates on the average of the Selberg trace
formula, which explains the increasing difficulty of this technique as we approach the optimal gap
$1/4$.

\subsection{Random covers and their spectral gap: from resonances to $L^2$-eigenvalues.}

The second author's interest for random covers was initially motivated by the spectral theory of
infinite volume surfaces. Let us assume that $X=\Gamma\backslash \H$ is a convex co-compact surface
i.e. $\Gamma$ is a finitely generated Fuchsian group of the second kind, without elliptic or
parabolic elements.  In that case $X$ has infinite volume, and is isometric to a compact surface
with geodesic boundary on which funnels have been glued (the infinite volume part).

\begin{figure}[htbp]
  \centering
  \includegraphics[scale=0.2]{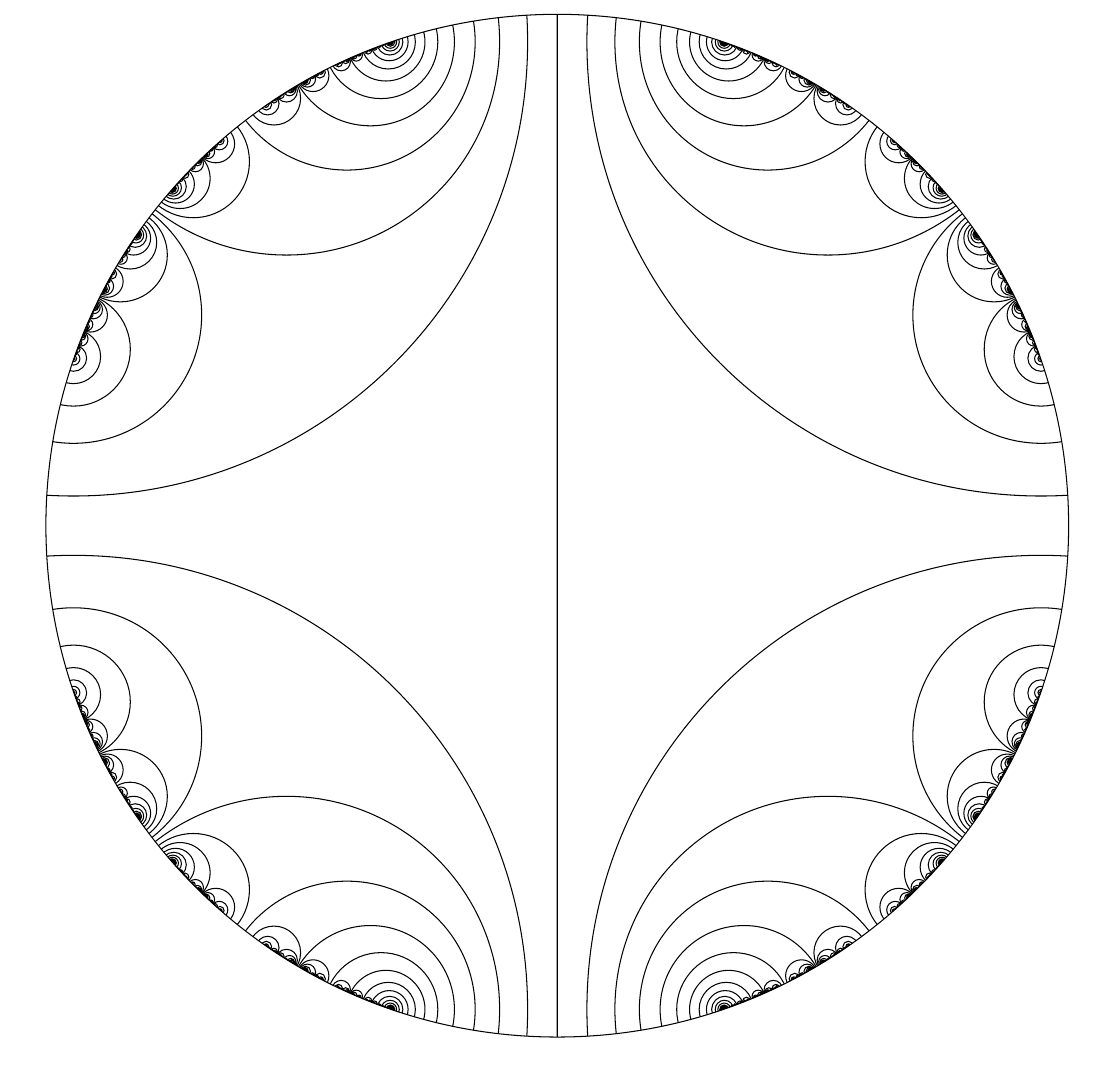}
  \includegraphics[scale=0.3]{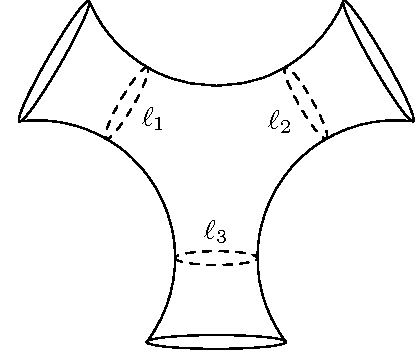}
  \caption{Limit set of Schottky group $\Gamma$ and pair of pants with funnels as quotient surface $X$.}
  \label{fig:schottky}
\end{figure}

In this situation, the $L^2$ spectrum of $\Delta_X$ has been described by Patterson and Sullivan and consists of absolutly continuous spectrum $[1/4,+\infty)$ and finitely many eigenvalues below $1/4$. See for example \cite{Borthwick} for an introduction to the subject and references herein. The bottom of the spectrum
is given by $\delta(1-\delta)$ if $\delta>1/2$, where $\delta$ is the Hausdorff dimension of the limit set of $\Gamma$, and is $1/4$ otherwise. 

If $s\in \C$ is Selbeg's spectral parameter, the resolvent 
$$R(s):=(\Delta_X-s(1-s))^{-1}:L^2(X)\rightarrow L^2(X),$$ as a function of $s$, is a holomorphic family of bounded operators on the
right half-plane $$\{ \Re(s)>1/2\},$$ except at possibly finitely many points on the real line, corresponding to the $L^2$- eigenvalues below $1/4$, where it is a finitely meromorphic operator valued map. By results of Mazzeo--Melrose \cite{mazzeomelrose} and Guillop\'e--Zworski \cite{GZ1}, one can show that the resolvent 
$$R(s):C_0^\infty(X)\rightarrow C^\infty(X),$$ has a meromorphic continuation to $\C$, and poles are called {\it resonances}, we will denote this set (with multiplicities) by $\mathcal{R}_X$. It is a known fact that resonances do coincides with zeros of the Selberg zeta function $Z_\Gamma(s)$, see \cite[Chapter 10]{Borthwick}, and resonances are actually intrinsic quantities which depend only on the length spectrum of $X$. The mathematical study of resonances in this hyperbolic setting is just an example of the broader picture of scattering theory and we recommend the survey of Zworski \cite{Zworski_survey} for a large overview of results and open conjectures.

At $s=\delta$, we always have a resonance and it is also known \cite{Naud1} that there exists $\epsilon(X)>0$ such that $\mathcal{R}_X\setminus \{\delta\}\subset \{ \Re(s)<\delta-\epsilon\}$. This leads to a natural definition of {\it spectral gap in the infinite volume setting} which is the size of the largest  possible $\epsilon(X)$. One can also relax this definition
and define also {\it the essential spectral gap} by 
$$\varepsilon_0(X):=\inf \{\sigma\leq \delta\ :\ \mathcal{R}_X\cap\{ \Re(s)\geq \sigma\}\ \mathrm{is\ finite} \}.$$
Bourgain and Dyatlov have shown in \cite{BD1} that there exists indeed $\varepsilon_0(\delta)<\delta$, which depends {\it only on $\delta$}, such that 
$\mathcal{R}_X\cap\{ \Re(s)\geq \sigma\}$ is finite. It was conjectured in \cite{JN1, JN2}, that one should have always $\varepsilon_0(X)=\varepsilon_0(\delta)=\frac{\delta}{2}$.

Motivated by the works of Gamburd \cite{Gamburd1}, Bourgain--Gamburd--Sarnak \cite{BGS} and Oh--Winter \cite{OW} on infinite index congruence-subgroups of $\mathrm{PSL}_2(\Z)$ (and a generalized Selberg conjecture in infinite volume), Michael Magee and the second author started investigating the behaviour of resonances in large degree covers in a random setting. One can define a notion of uniformly random covers of degree $n$ as follows. Let $X=\backslash \H$ be a fixed based surface, where $\Gamma$ is a finitely generated Fuchsian group, and let $\phi_n:\Gamma \rightarrow S_n$
be a homomorphism from $\Gamma$ to the Symmetric group of permutations of order $n$. Define $X_n$ to be quotient of
$\H\times\{1,\ldots,n\}$ by the action of $\Gamma$ given by $$\gamma.(z,j):=(\gamma(z),\phi_n(\gamma)(j)),$$
for all $\gamma \in \Gamma$. There is a natural covering map $X_n\rightarrow X$, and the surface $X_n$ is isometric to finite union of hyperbolic surfaces
$$X_n=\bigsqcup_k \Gamma_k \backslash \H,$$
where all subgroups $\Gamma_k\subset \Gamma$ have finite index and $\sum_k [\Gamma: \Gamma_k]=n$. Since the space $\mathrm{Hom}(\Gamma,S_n)$ is finite, we
endow it with the {\it uniform} probability measure and obtain a well defined notion of random cover. Note that the cover may not be connected, but this occur with a probability that
usually tends to $0$ as $n\rightarrow +\infty$ (for example as a consequence of \cite{broder1987} for the free group). When $X$ is a non-compact surface, $\Gamma$ is a {\it free group}, which makes it possible to use the combinatorials tools developped for random graphs \cite{broder1987, Puder1}. One of the key facts is the asymptotics for the number of fixed points:
given $\gamma \in \Gamma$ which is primitive word (i.e. $\gamma$ is not a power $\gamma=\gamma_0^k$ for some $k\geq 2$), we have as $n\rightarrow +\infty$, the expectation asymptotics
\begin{equation}
\label{Fix_average}
\mathbb{E}_n(\mathrm{Fix}(\phi_n(\gamma)))=1+ O_\gamma\left(\frac{1}{n}\right),
\end{equation}
where $\mathrm{Fix}(\sigma)$ is the number of fixed points of the permutation $\sigma \in S_n$. We point out that for uniformly distributed random permutations $\alpha\in S_n$, it is elementary to check that $\mathbb{E}(\mathrm{Fix}(\alpha))=1$. 

The constant implied in the bound (\ref{Fix_average}) depends on the word length $\vert \gamma\vert$ of $\gamma$. If the word length has moderate growth with respect to $n$, i.e. $\vert \gamma \vert\leq C_1\log(n)$, then it is possible to find a uniform constant $C_2$ such that 
$$\vert \mathbb{E}_n(\mathrm{Fix}(\phi_n(\gamma)))-1\vert \leq C_2\log(n)/n.$$
Combining these combinatorial tools with Hilbert--Schmidt estimates for transfer operators from \cite{JN2}, Michael Magee and the second author \cite{MN1} were able to show the following fact. Set 
$$M(\sigma,T):=\{ s \in \C\ :\ \sigma\leq \Re(s)\leq \delta\ \mathrm{and}\ \vert\Im(s)\vert\leq T\}.$$

\begin{theorem}
\label{gap1}
Fix $T>0$ and pick $\frac{3}{4}\delta<\sigma<\delta$, then with high probability as $n\rightarrow +\infty$, we have
$$\mathcal{R}_{X_n}\cap M(\sigma, T)=\mathcal{R}_{X}\cap M(\sigma, T).$$
\end{theorem}
In other words, there are no new resonances in the box $M(\sigma,T)$ with a probability which tends to $1$ as $n$ goes to infinity. Note that resonances or eigenvalues of the Laplacian on $X$ are always a subset (with multiplicities) of $\mathcal{R}_{X_n}$. We point out that if formally $\delta=1$, and taking in account the original spectral parameter $s(1-s)$, then
the gap obtained is the same as Selberg's $3/16$. The proof we gave in \cite{MN1} can essentially be considered as an infinite dimensional version of the {\it trace method} used in
\cite{broder1987, Puder1} for random permutation graphs. 

Let us now go back to finite volume or compact surfaces. We recall that all the previous uniform
bounds for spectral gaps were proved by applying Cheeger's inequality, and that this approach is
limited and cannot get close to the optimal value $1/4$ for the spectral gap. It became clear at this
point that to address the problem of compact surfaces, a surface group analog of \cite{Puder1} was
needed. This was achieved by Michael Magee and Doron Puder first in \cite{MP1}, and refined for the
problem at hand in \cite{magee2022} to obtain the following result.

\begin{theorem} Assume that $X$ is a compact surface.
\label{gap2}
For any $0<\lambda<\frac{3}{16}$, with high probability as $n\rightarrow +\infty$, we have
$$\mathrm{Sp}(X_n)\cap [0,\lambda]=\mathrm{Sp}(X)\cap [0,\lambda],$$
where $\mathrm{Sp}(X)$ is the (discrete) $L^2$-spectrum of the Laplacian on $X$.
\end{theorem}
The proof is based on the combination of Selberg's trace formula (with suitable test functions) and the expectation asymptotics for the number of fixed points as above, in perfect analogy with the trace method for graphs. The hard part is to prove the analog of (\ref{Fix_average}) for the combinatorics of surface groups.

\subsection{Weil--Petersson random surfaces and spectral gap.}

Shortly after \cref{gap2} was written for random covers of a fixed hyperbolic surface, Wu and Xue
\cite{wu2022} and independently Lipnowski and Wright \cite{lipnowski2021} proved an analogous
$3/16$-result in the large genus regime for the Weil-Petersson model (see also \cite{hide2022}).

\begin{theorem}
  \label{gap316wp}
  For any $0 < \lambda < \frac{3}{16}$, we have $\lambda_1(X) > \lambda$ with high probability as
  $g \rightarrow + \infty$.
\end{theorem}

In both articles, the result is obtained by a similar trace method using the Selberg trace formula,
which is also shared by the random cover article \cite{magee2022} as it is a natural way to tackle
this problem. The authors use a dilated test function $h_L (r) = h(r/L)$ for a large value of $L$,
identify a contribution corresponding to the trivial eigenvalue, and prove the rest is small enough
to reach the desired gap. The non-optimal value $3/16$ instead of $1/4$ appears as a consequence of
the level of precision all the estimates are carried on, $1/g$ or $1/n$ respectively.

This is where the similarities end, and the proof becomes model-specific.  In the case of the
Weil--Petersson model, instead of counting fixed points of random permutations as in
\cref{Fix_average}, one needs to average quantities of the form
$\sum_{\gamma \in \mathcal{P}_X} F(\ell(\gamma))$, i.e. sums over all primitive closed
geodesics. This should be reminiscent of the geometric functions studied by Mirzakhani, which, as we
saw in \cref{sec:wolp-mirz-magic}, can be integrated. Unfortunately, there is a significant sticking
point here: Mirzakhani's integration formula only applies to \emph{simple} closed geodesics (with no
self-intersections), whilst \emph{all} closed geodesics appear in the Selberg trace formula.

At this point, both teams make a stricking observation. If one is to restrict the Selberg trace
formula to only simple closed geodesics, and integrate it using Mirzakhani's techniques, they
recover the constribution of the trivial eigenvalue $\lambda_0=0$, up to a small error. More
precisely, 
\begin{equation}
  \label{e:simple_match}
  \mathbb{E}_g^{\mathrm{WP}} \Bigg[ \sum_{\substack{\gamma \in \mathcal{P}_X \\ \text{simple}}}
  \frac{\ell h_L(\ell(\gamma))}{2 \sinh (\ell(\gamma)/2)}\Bigg]
  = 2 \int_0^L h_L(\ell) \sinh \Big( \frac \ell 2 \Big) \d \ell
  + \mathcal{O}_\epsilon \Big(\frac{e^{(1+\epsilon)L/2}}{g} \Big)
\end{equation}
where the first term is very close to $h_L(i/2) = 2 \int_0^L h_L(\ell) \cosh (\ell/2) \d \ell$ for
large values of $L$.  This observation is parallel to \cref{Fix_average} in the case of random
covers, but comes from completely distinct considerations.

The last (and most substantial) step then becomes to justify that the other contributions are small.
Recall that the Selberg trace formula here contains all closed geodesics of length $\leq L$, for a
large value of $L$. By Huber's counting result, we know that there are exponentially many such
geodesics, which means that the number of terms appearing in the sum is a priori way too large to
say anything. We need to obtain more information on this exponential growth, in particular for
random hyperbolic surfaces and the specific $L$ we choose (here, $L = 4 \log g$).

\begin{figure}[htbp]
  \centering
  \begin{subfigure}[b]{0.45\textwidth}
    \centering
    \includegraphics[scale=0.4]{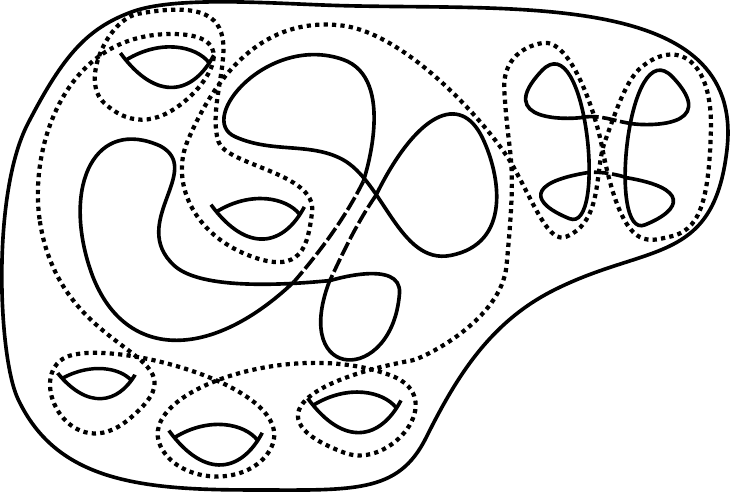}
    \caption{Filling a large subsurface.}
    \label{fig:longA}
  \end{subfigure}%
  \begin{subfigure}[b]{0.45\textwidth}
    \centering
    \includegraphics[scale=0.4]{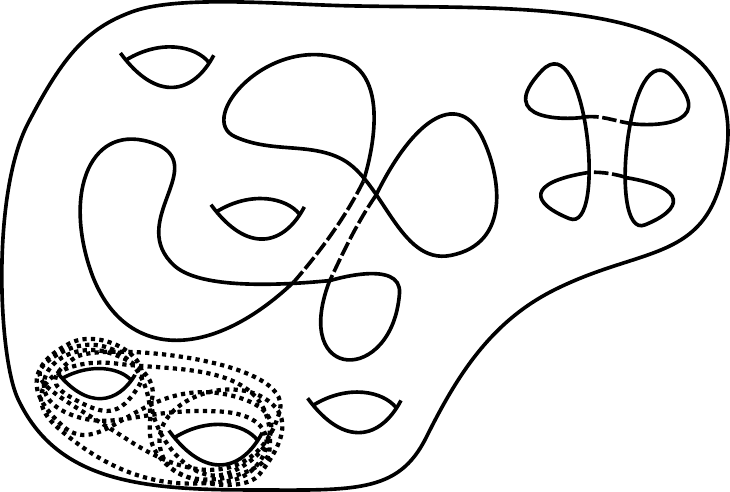}
    \caption{Filling a pair of pants.}
    \label{fig:longB}
  \end{subfigure}%
  \caption{Two long closed geodesics on a large-genus surface. }
  \label{fig:long}
\end{figure}

Long geodesics on a compact hyperbolic surface can do different things, represented in Figure
\ref{fig:long}. Some might explore a very large chunck of the surface (as in Figure
\ref{fig:longA}), whilst other might be constrained to a small part, but wind around many times (as
in Figure \ref{fig:longB}). It is therefore interesting to introduce the \emph{surface filled by a
  closed geodesic}, as the minimal surface with geodesic boundary containing it. Mirzakhani's
estimates allow to show quite straightforwardly that geodesics for which the Euler characteristic of
the filled surface is very negative (such as in Figure \ref{fig:longA}) correspond to very small
contributions on average and can be discarded.

The difficulty therefore comes from geodesics winding down many times in a small section of the
surface, as in Figure \ref{fig:longB}. These occur in particular when the surface contains
\emph{tangles}, which are embedded surfaces of Euler characteristic $-1$ with a short boundary
length. This notion has been introduced by the first author and Joe Thomas in \cite{monkTF},
adapting a similar concept from the works of Friedman and Bordenave on the Alon
conjecture. \emph{Tangled surfaces} (surfaces containing tangles) have too many closed geodesics,
which causes the trace method to fail. They occur with small probability, but yield large
contributions on average.  This challenge is overcome in two very different ways in \cite{wu2022}
and \cite{lipnowski2021}, which each brings insight on the geometry of random Weil--Petersson
surfaces and we shall now detail.

\begin{itemize}
\item On the one hand, Wu and Xue prove a refined upper bound on the number of closed geodesics
  filling a given surface. Let $Y$ be a surface with geodesic boundary. Then, Wu and Xue prove that
  the number of closed geodesics of length $\leq L$ filling $Y$ is bounded above by
  $e^{L - (1-\epsilon) \ell(\partial Y)/2}$. In other words, the Huber bound can be improved with a
  correction decaying exponentially with the length of the boundary of $Y$. This refinement is
  exactly what is needed to reach the level of precision required for the $3/16$ result.
\item On the other hand, Lipnowski and Wright prove a polynomial bound on the number of closed
  geodesics of length $\leq L$ contained in a tangle-free surface. The shift from exponential to
  polynomial bound is enough to drastically reduce the difficulty of the argument, with one major
  trade-off. They now need to integrate on a subset of the moduli space, here the \emph{thick part}
  $\mathcal{M}_g^{> \kappa}$ (hyperbolic surfaces with an injectivity radius $> \kappa$). We
  have seen in \cref{sec:wolp-mirz-magic} that the theory of integration on the moduli space is
  challenging, and the types of functions we can integrate are strongly prescribed. Luckily, one can
  rewrite the indicator function of the thick part by a simple inclusion-exclusion, as first
  performed in \cite{mirzakhani2013}:
  \begin{equation}
    \label{eq:inc_excl}
    \mathds{1}_{\mathcal{M}_g^{> \kappa}}(X)
    = 1-\sum_{q=1}^\infty \frac{(-1)^{q+1}}{2^q q!}
    \sum_{(\gamma_1, \ldots, \gamma_q)}
    \prod_{j=1}^q \mathds{1}_{[0,\kappa]}(\ell(\gamma_j)).
  \end{equation}
  The sum above runs over all families of distinct primitive closed geodesics $(\gamma_1, \ldots,
  \gamma_q)$. Conveniently, due to the collar lemma, if $\kappa$ is small enough, then all of the
  geodesics contributing to this sum are simple and disjoint, and hence the formula above can be
  viewed as a rewriting of $\mathds{1}_{\mathcal{M}_g^{> \kappa}}$ as a combination of geometric
  functions, which can then be integrated using Mirzakhani's integration formula.
\end{itemize}

This $\frac{3}{16}$ result is a substantial improvement on the initial bound $0.002$ proven by
Mirzakhani using the Cheeger inequality. We know from Friedman and Bordenave's proofs of the Alon
conjecture that trace methods have the potential to reach the optimal spectral gap $\frac
14$. Nalini Anantharaman and the first author achieved this in 2025:

\begin{theorem}
  \label{thm:WP}
  For any $0 < \lambda < \frac 14$, we have $\lambda_1(X) > \lambda$ with high probability as
  $g \rightarrow + \infty$.
\end{theorem}

The proof of this result appears as a two-part paper \cite{anantharaman2023,anantharaman2025}
totalling 263 pages, which is summarised in an exposition article \cite{anantharaman2024}. This jump
in complexity is due to the fact that the number $\frac{3}{16}$ appears as a precision barrier for
multiple different reasons, which all need to be tackled to improve on it, and ultimately reach the
optimal gap~$\frac 14$.  Let us describe a few aspects of the approach.

The spectral gap $\frac{3}{16}$ is obtained by computing the average of the Selberg trace formula at
the leading order, i.e. with error terms decaying in $1/g$. A little computation allows to see that,
in order to obtain a better gap using the trace method, it is necessary to increase the precision in
the calculations (see \cite[Section 3.4]{anantharaman2023}). More precisely, the best gap we can
hope for if computing all terms with errors decaying in $1/g^N$ is $\frac 14 -
\frac{1}{4(N+1)^2}$. As a consequence, the only way to reach $\frac 14$ via the trace method
requires asymptotic expansions in powers of $1/g$ with arbitrary level of precision. Certain
expansions of Weil--Petersson volumes were derived by Mirzakhani and Zograf in
\cite{mirzakhani2015}; more information on their coefficients were obtained in
\cite{anantharaman2022}.

As we compute expectations with higher and higher level of precision, more and more complex closed
geodesics start to appear. For instance, the figure-eight (geodesics with exactly one
self-intersection) contribute to the term $1/g$. This is arguably the biggest challenge arising from
the need for an expansion, as the only tool available, Mirzakhani's integration formula, is
fundamentally limited to the study of simple geodesics. Anantharaman and the first author defined in
\cite{anantharaman2023} an equivalence relation on closed curves called \emph{local topological
  type} $\mathbf{T}$, which allows to group closed geodesics of similar topology. An explicit
formula for the \emph{volume function} associated to a given local topological type, i.e. the function
$V_g^{\mathbf{T}} : \R_{>0} \rightarrow \R$ such that
\begin{equation}
  \label{eq:type}
  \mathbb{E}_g^{\mathrm{WP}} \Big[ \sum_{\substack{\gamma \in \mathcal{P}_X \\ \text{of type }
      \mathbf{T}}}
  F(\ell(\gamma)) \Big] = \int_0^\infty F(\ell) V_g^{\mathbf{T}}(\ell) \d \ell,
\end{equation}
is derived in \cite{anantharaman2023}. A full asymptotic expansion in powers of $1/g$ is
obtained. The formula contains a somehow abstract integration on a Teichm\"uller space (allowing to
list all possible geometries for the closed geodesic). It is made explicit and fully computable by
the introduction of novel sets of coordinates on Teichm\"uller spaces, in which the Weil--Petersson
form admits a surprisingly simple expression, in \cite{anantharaman2025}.

The miracle correspondance between the contribution of the trivial eigenvalue and that of simple
closed geodesics, observed in \cref{e:simple_match}, only occurs at the leading order. It is quite
unclear how one would proceed to prove such a precise cancellation at arbitrary high orders of
precision. For this reason, a different, more robust approach is used to cancel the trivial
eigenvalue $0$ in \cite{anantharaman2023,anantharaman2025}. The test function $h_L$ is replaced by a
new test function which has a Fourier transform vanishing at $r_0=i/2$. This is done by applying
a power of the differential operator $\frac 14 - \partial^2$ to the test function $h_L$, which, by
Fourier analysis, corresponds to multiplying the Fourier transform by a power of $\frac 14 +
r^2$. The objective then becomes to exhibit non-trivial cancellations in the geometric term of the
trace formula. The argument is detailed in \cite[Section 3.4]{anantharaman2023} and leads to the
introduction of a fundamental assumption:

\begin{definition}
  A locally integrable function $f : \R \rightarrow \R$ is called a \emph{Friedman--Ramanujan
    function} if there exists a polynomial $p$ such that $f(\ell) = p(\ell) e^\ell +
  \mathcal{O}((\ell+1)^c e^{\ell/2})$.
\end{definition}

This expression can be compared to the prime number theorem \cref{eq:prime}: the main term
$p(\ell)e^\ell$ appears as coming from the trivial eigenvalue $0$, and the gap between the exponents
$e^\ell$ and $e^{\ell/2}$ corresponds to the gap between $0$ and $\frac 14$. In effect, terms of the
form $p(\ell) e^\ell$ are exactly cancelled in the Selberg trace formula by a simple integration by
parts, using the differential operator $\frac 14 - \partial^2$.
The main result of \cite{anantharaman2025} can then be reformulated as follows.
\begin{theorem}
  For any non-simple local topological type $\mathbf{T}$, any order of precision $N$, the $N$-th
  order approximation of the volume function $V_g^{\mathbf{T}}$ is a Friedman--Ramanujan function in
  the weak sense.
\end{theorem}
The full result contains more information on the uniformity of the degrees, constants and
coefficients of the Friedman--Ramanujan functions, which is required to obtain a spectral gap. The
weak sense refers to a weaker, integrated, definition of Friedman--Ramanujan functions.

The last significant challenge tackled in the proof of \cref{thm:WP} comes from the necessity to
remove tangles, as mentioned above in Lipnowski and Wright's approach. There it is solved following
the lines of Friedman's work. The basic idea is to rewrite the indicator function of the set of
tangle-free surfaces as a sum of geometric functions, similarly to \cref{eq:inc_excl}. However, past
the precision $\frac{3}{16}$, one needs to include tangles (surfaces of Euler characteristic $-1$
with a short boundary) in this inclusion-exclusion. Contrarily to short closed geodesic, these
tangles can intersect, and one hence a priori needs to list all possible combinations of
tangles. This issue is solved by proving a \emph{Moebius formula} in \cite{anantharaman2024}, which
achieves exactly the goal stated above, whilst avoiding tedious topological enumerations.

\section{The strong convergence breakthrough.}

Apart from the trace method, a recent breakthrough paper of Bordenave-Collins \cite{bordenave2019} had caught the attention of Michael Magee and the second author.
The main result of \cite{bordenave2019} is about {\it strong convergence} for sums of random
permutation matrices, it implies in particular Alon's conjecture for random lifts of regular
graphs. Strong convergence is an alternative approach to the trace method and is in some situations
more effective.

\subsection{Strong convergence.}

Let us explain first what strong convergence is about. Let $\Gamma$ be a {\it finitely generated discrete group}, and let $\C[\Gamma]$ be the group ring of $\Gamma$ i.e. the ring of formal finite sums $\sum_{\gamma \in \Gamma} a_\gamma$, where $a_\gamma\in \C$. Let $\rho_n:S_n\rightarrow \mathrm{End}(L^2(\{1,\ldots,n\})$
be the standard representation of the symmetric group $S_n$ defined for $\sigma \in S_n$ by 
$$\rho_n(\sigma)(f)(j):=f(\sigma^{-1}(j)),$$
where $f\in L^2(\{1,\ldots,n\})\simeq \C^n$. Let $V_n^0\subset L^2(\{1,\ldots,n\})$ be the subspace of vectors which are orthogonal to constants functions i.e.
$$V_n^0=\Big\{ f\in L^2(\{1,\ldots,n\})\ :\ \sum_j f(j)=0 \Big\}. $$
Then the restriction of $\rho_n$ to $V_n^0$ is an irreducible representation of $S_n$, which will refer as $\rho_n^0$ in the following.

A sequence of permutation representations $(\phi_n)$ with $\phi_n:\Gamma \rightarrow S_n$ is said to converge strongly if, for any element $z\in \C[\Gamma]$, 
we have 
$$\lim_{n\rightarrow +\infty} \Vert \rho_n^0(\phi_n(z))\Vert_{V_n^0}=\Vert \lambda_\Gamma(z) \Vert_{\ell^2(\Gamma)},$$
where the norms involved are operator norms, and
$\rho_n^0(\phi_n(z))=\sum_{\gamma}a_\gamma \rho_n^0(\phi_n(\gamma))$, if $z=\sum_\gamma a_\gamma$,
and $\lambda_\Gamma(z)$ is the left regular representation acting on $\ell^2(\Gamma)$. A more
pedestrian way to say it: on the orthogonal of constant vectors, (the norm of ) the sum of
permutation matrices converges to the norm of the left regular representation of~$z$. In
\cite{bordenave2019}, they showed that if {\it $\Gamma$ is the free group} and
$\phi_n:\Gamma\rightarrow S_n$ is sampled according to the uniform probability measure on
$\mathrm{Hom}(\Gamma,S_n)$, then strong convergence occurs with probability tending to $1$ as $n$
goes to infinity.

It is important to point out that strong convergence actually implies convergence when tensoring with finite dimensional matrices, a trick known as ``matrix amplification'', which can be deduced from strong convergence, for more details see \cite{magee_survey}, section $3$ and references herein. 

\subsection{Application to finite volume covers and compact surfaces with near optimal spectral
  gap.}

In \cite{MN2}, Michael Magee and the second author managed to combine the non-backtracking approach from \cite{bordenave2020} and the transfer operator techniques from \cite{MN1} to prove the following fact.

\begin{theorem}
\label{gap2}
Let $X=\Gamma \backslash \H$ be a convex co-compact base surface as above, where $\Gamma$ is a \underline{free group}.
Fix $T>0$ and pick any $\sigma$ with $\frac{\delta}{2}<\sigma<\delta$, then with high probability as $n\rightarrow +\infty$, we have
$$\mathcal{R}_{X_n}\cap M(\sigma, T)=\mathcal{R}_{X}\cap M(\sigma, T).$$
\end{theorem}
The $\delta/2$-bound obtained above is now the {\it infinite volume} analog of the $1/4$-gap stated above for compact and finite volume surfaces, and proves part of the conjecture from \cite{JN2}. The proof given at the time in \cite{MN2} is fairly technical and has since been made more conceptual in \cite{MCN}, where we combine directly the strong convergence result from Bordenave-Collins and Haagerup inequalities to estimate the norm of the limit operator. The fact that strong convergence methods could lead directly to near-optimal result in the infinite volume setting was a clear indicator that perhaps the finite volume case was reachable. This is exactly what Hide and Magee achieved in \cite{hide2021}.
\begin{theorem}
\label{gap3}
Let $X=\Gamma \backslash \H$ be a finite volume base surface, in particular $\Gamma$ is a \underline{free group}.
For all $\epsilon>0$,  with high probability as $n\rightarrow +\infty$, we have
$$\mathrm{Sp}(X_n)\cap [0,1/4-\epsilon]=\mathrm{Sp}(X)\cap [0,1/4-\epsilon],$$
where $\mathrm{Sp}(X)$ is the $L^2$-spectrum of the Laplacian on $X$.
\end{theorem}
For example, if the base surface $X$ has Euler characteristic $\chi(X)=-1$, then by a bound of Otal and Rosas \cite{otal2009}, $0$ is the only eigenvalue below $1/4$ and thus for large $n$,
with high probability we have
$$\mathrm{Sp}(X_n)\cap [0,1/4-\epsilon]=\{0\}.$$
In the same paper, Hide and Magee explain how the arguments from \cite{buser1988} allow to compactify a sequence of finite volume surfaces with near optimal spectral, while controling the genus and without loosing to much on $\lambda_1$.  Therefore, they obtained as a by product the first proof of existence of a sequence of compact surfaces $(X_k)$ with genus going to infinity and $\lim_{k\rightarrow \infty}\lambda_1(X_k)= 1/4$. The proof in \cite{hide2021} is based on an explicit parametrix for the resolvent $R_{X_n}(s)$ of the 
Laplacian on the covers $X_n$. This parametrix involves a (finite) summation on the group $\Gamma$ and is written as a sum of permutation matrices tensored by compact operators. Applying directly the strong convergence result allows to relate (with high probability) norm estimates of the parametrix to norm estimates of the resolvent on the universal cover $\H$, for which the spectrum is known to be $[1/4,+\infty)$. Michael Magee has recently managed to show \cite{magee_survey}, that in {\it finite volume locally symmetric covers}, a proof of the near-optimal spectral gap can be achieved by pure representation theoretic arguments, thus eleminating all non-necessary specific analysis in the proofs.

We conclude this paragraph by pointing out that, at the time of writing of the present survey, the optimal spectral gap result is not known for the Brooks-Makeover model. There are two main reasons for it. First Brooks-makeover surfaces arise from picking gluing combinatorics out of finite index subgroups of free products of cyclic groups (one of them has order $3$). Therefore Bordenave-Collins result cannot be applied directly. Secondly, the uniform probability mesure on the space of finite index subgroup of this free product differs slightly from the Brooks-Makeover probability measure. It is however likely that by using the new polynomial approach to strong convergence \cite{chen2025} one could reach the $1/4$-gap.

\subsection{Strong convergence for surface groups.}
In view of the powerful applications derived from strong convergence for sums of permutations matrices, a tantalizing question is can one prove strong convergence for surface groups? So far strong convergence was only known for free groups. A first partial answer is given by Louder and Magee in \cite{louder2022} where they are able to show that for any abstract surface group $\Gamma$, there exists a sequence
of permutation representations $(\phi_n):\Gamma\rightarrow S_n$, such that the corresponding sequence of permutation matrices (acting on $V_n^0$) strongly converges to the left regular representation on $\ell^2(\Gamma)$. Their proof is based on "perturbative arguments" and relies on the free case and Bordenave-Collins result. The novelty is that one can start from any base compact surface and obtain a sequence of covers with near optimal relative spectral gap. Starting for example from $X$ being Bolza's surface, which is arithmetic and known to have $\lambda_1\simeq 3.8388$, they obtain the existence of a sequence of {\it compact arithmetic surfaces} with genus going to infinity such that $\lim_{k\rightarrow \infty}\lambda_1(X_k)= 1/4$.

\subsection{The polynomial method.}
\label{s:polynomial}

Recently, an alternative proof of strong convergence (for the free group) was proposed by Chen, Garza-Vargas, Tropp and Van Handel \cite{chen2025}, showing that one could use soft arguments to ``bootstrap'' the moments asymptotics (i.e. results known as asymptotic freeness in the litterature). The key input is the fact that for many models (random permutations matrices, random Haar unitary matrices) the moments asymptotics are rational functions of $1/N$, where $N$ is the matrix size. The technique used is then polynomial approximation theory and basic Fourier analysis. 

This new approach, known as ``the polynomial method'' has opened up new perspectives on strong convergence. A preprint from 2025 by Magee, Puder and Van Handel \cite{magee2025} settles the case of surface groups: they managed to produce a proof of strong convergence by the polynomial method. There are key new ideas in this paper since 
the rationality of moments functions fails for surface groups: they work instead in Gevrey regularity class, but still are able to achieve the tour de force of proving strong convergence.
We can therefore state the following main consequence.
\begin{theorem}
\label{gap4}
Let $X=\Gamma \backslash \H$ be a compact hyperbolic surface.
For all $\epsilon>0$,  with high probability as $n\rightarrow +\infty$, we have
$$\mathrm{Sp}(X_n)\cap [0,1/4-\epsilon]=\mathrm{Sp}(X)\cap [0,1/4-\epsilon],$$
where $\mathrm{Sp}(X)$ is the $L^2$-spectrum of the Laplacian on $X$.
\end{theorem}
Let us mention that the strong combinatorial result from \cite{magee2025} has the potential for applications far beyond hyperbolic surfaces. We just mention here a recent work
of Moy, Hide and the second author \cite{HMN}, where we obtain the analog result for random covers of (variable) negatively curved surface covers, by combining \cite{magee2025}
and Heat kernel techniques.
More recently, Hide, Macera and Thomas successfully implemented ideas of the polynomial method to the
Weil--Petersson model \cite{hide2025}, and hence obtained an alternative proof with a spectral gap with a polynomial rate
$\frac 14 - \frac{1}{g^c}$ for a $c>0$.

\section*{Acknowledgments.}
We both thank Bram Petri for his careful reading and comments. 
This research was funded by the EPSRC grant EP/W007010/1 and the Royal Society Dorothy Hodgkin
Fellowship for the first author.

\bibliographystyle{plain}

\end{document}